\documentclass[12pt]{article}[12pt]
\usepackage[latin1]{inputenc}
\usepackage[english]{babel}
\usepackage{amsmath,amssymb}
\newcommand{\bbR} {{\mathbb{R}}}   
   
\newcommand{\bbC} {{\mathbb{C}}}   
\newcommand{\bbK} {{\mathbb{K}}}   
\newcommand{\bbP} {{\mathbb{P}}}

\newtheorem{Theorem}{Theorem}[section]

\newtheorem{Lemma}[Theorem]{Lemma}
\newtheorem{Corollary}[Theorem]{Corollary}

\newtheorem{Example}[Theorem]{Example}

\def\ftoday{le\space\number\day \space\ifcase\month\or
  janvier\or f\'evrier\or mars\or avril\or mai\or juin\or
  juillet\or ao\^ut\or septembre\or octobre\or novembre\or d\'ecembre\fi
  \space\number\year}

\newcommand{\cqfd}{\hfill$\Box$}

\begin{document}

\title{Transcendental versions in $\bbC^n$ of the Nagata conjecture}
\author{St\'ephanie Nivoche\thanks{Universit\'e C\^ote d'Azur, Laboratoire J.A. Dieudonn\'e, 
CNRS UMR 7351, Universit\'e de Nice Sophia Antipolis, Parc Valrose, 06108 Nice cedex 2, France. 
E-mail : nivoche@unice.fr. Research of S. Nivoche was supported by ANR grant ``ANR-14-CE34-0002-01'' and the FWF grant I1776 for the international cooperation project ``Dynamics and CR geometry''.}}
\maketitle
\rmfamily

\begin{abstract}
The Nagata Conjecture is one of the most intriguing open problems in the area of curves in the plane. It is easily stated. 
Namely, it predicts that the smallest degree $d$ of a plane curve passing through $r \ge 10$ general points in the projective plane $\bbP^2$ with multiplicities at least $l$ at every point, satisfies the inequality
$d > \sqrt{r} \cdot l$. This conjecture has been proven by M. Nagata in $1959$, if $r$ is a perfect square greater than $9$. 
Up to now, it remains open for every non-square $r \ge 10$, after more than a half century of attention by many researchers.\\
In this paper, we formulate new transcendental versions of this conjecture coming from pluripotential theory and which are equivalent to a version in $\bbC^n$ of the Nagata Conjecture.
\end{abstract}



\newpage

\section{Introduction}

\subsection{History and known results}

In $1900$, at the international Congress of Mathematicians in Paris, D. Hilbert posed twenty-three problems. 
The fourteenth one may be formulated as follows: 
{\it Let $\bbK$ be a field and $x_1 , \ldots , x_n$ algebraically independent elements 
over $K$. Let $L$ be a subfield of $\bbK(x_1 , . . . , x_n)$ containing $\bbK$. Is the ring $\bbK [x_1 , . . . , x_n ] \cap L$ finitely generated over $\bbK$ ?}\\
Hilbert conjectured that all such algebras are finitely generated over $\bbK$. \\
Contributions to the fourteenth problem are obtained confirming Hilbert's conjecture in special cases and for certain classes of rings. 
In $1953$, a significant contribution was made by 0. Zariski, who generalized
the fourteenth problem in the following way :\\
\indent  {\it Problem of Zariski. Let $\bbK$ be a field and $\bbK [a_1 ,\ldots , a_n ]$ an affine
normal domain (i.e. a finitely generated integrally closed domain over
$\bbK$). Let $L$ be a subfield of $\bbK (a_1 , \ldots, a_n )$ containing $\bbK$. Is the ring $\bbK
[a_1 ,\ldots, a_n ]\cap  L$ finitely generated over $\bbK$ ?}\\
He answered the question in the affirmative when trans.deg$_{\bbK} L \le 2$. Later, in $1957$, D. Rees 
gave a counter example to the problem of Zariski when trans.deg$_{\bbK} L = 3$. Finally Masayoshi Nagata in \cite{Nag2}
gave a counter example to the original fourteenth problem itself. This counter example is in
the case of trans.deg$_{\bbK} L = 13$. \\

In $1959$, M. Nagata \cite{Nag1} gave another counter example (a suitably constructed ring of invariants for the action of a linear algebraic group) in the case of trans.deg$_{\bbK} L = 4$.\\
In this work, he finally formulated a conjecture, which governs the minimal degree required for a plane algebraic curve to pass through a collection of general points with prescribed multiplicities : 
{\it Suppose $p_1, \ldots, p_r$ are $r$ general points in $\bbP^2$ and that $m_1,\ldots, m_r$ are given positive integers. 
Then for $r > 9$, any curve $C$ in $\bbP^2$ that passes through each of the points $p_i$ with multiplicity $m_i$ must satisfy $\deg C > \frac{1}{\sqrt{r}} \cdot \sum_{i=1}^r m_i.$}  \\
One says that a property ${\cal P}$ holds for $r$ general points in $\bbP^2$ if there is a Zariski-open subset $W$ of $(\bbP^2)^r$ such that ${\cal P}$ holds for every set $(p_1,\ldots,p_r)$ of 
$r$ points in $W$.\\
As Nagata pointed out, it is enough to consider the uniform case. Thus this conjecture is equi\-valent to the following one, which it usually called :\\

{\it The Nagata Conjecture. Suppose $p_1, \ldots, p_r$ are $r$ general points in $\bbP^2$ and that $m$ is a given positive integer. 
Then for $r > 9$, any curve $C$ in $\bbP^2$ that passes through each of the points $p_i$ with multiplicity at least $m$ must satisfy $\mbox{deg}\, C > \sqrt{r} \cdot m.$}\\

The only case when this is known to hold true is when $r>9$ is a perfect square. This was proved by Nagata with technics of specializations of the points in the plane. 
More recently, without any extra condition connecting the degree, $m$ and $r$, G. Xu \cite{Xu} proved that $\deg C \ge \frac{\sqrt{r-1}}{r} \cdot \sum_{i=1}^r m_i$ and H. Tutaj-Gasinska \cite{TG} proved that 
$\deg C \ge \frac{1}{\sqrt{r+(1/12)}} \cdot \sum_{i=1}^r m_i$ (see also B. Harbourne \cite{Ha2} and B. Harbourne and J. Ro\'e \cite{HaRo}). \\

A more mo\-dern formulation of this conjecture is often given in terms of Seshadri constants, introduced by J.-P. Demailly (\cite{De5} and \cite{De6}) in the course of his work on Fujita's conjecture.
The Nagata Conjecture is generalized to other surfaces under the name of the Nagata-Biran conjecture (let $X$ be a smooth algebraic 
surface and $L$ be an ample line bundle on $X$ of degree $d$. The Nagata-Biran conjecture states that for sufficiently large $r$ the Seshadri constant satisfies
$\epsilon(p_1,\ldots,p_r;X,L) = \frac{d}{\sqrt{r}}$).\\
Nagata has also remarked that the condition $r > 9$ is necessary. The cases $r > 9$ and $r \le 9$ are distinguished by whether or not the anti-canonical bundle on the blowup of $\bbP^2$ at a 
collection of $r$ points is nef.\\
Fix $r$ general points $p_1, \ldots, p_r$ in $\bbP^2$ and a nonnegative integer $m$. Define $\delta(r,m)$ to be the least integer $d$ such that there is a curve of degree $d$ vanishing at each 
point $p_i$ with multiplicity at least $m$.
For $r \le 9$, applying methods of \cite{Nag3}, Harbourne \cite{Ha2} shows that $\delta(r,m)=\lceil c_r m \rceil$, where $c_r=1; \, 1; \, 3/2; \, 2; \, 2; \, 12/5; \, 21/8; \, 48/17$ and $3$ 
for $r=1,\ldots,9$ respectively (for any real number $c$, $\lfloor c \rfloor$ is the greatest integer less than or equal to $c$). 
When $r=1, 4$ and $9$, we can remark that the Nagata conjecture holds with a weak inequality instead of a strict one.\\

Iarrobino \cite{Iar} The generalized Nagata Conjecture in any dimension $n \ge 2$. Based on a conjecture of Fr\"{o}berg,  Iarrobino predicted that : {\it an hypersurface in $\bbP^n$ passing 
through $r$ generic 
points with multiplicity $m$ has a degree $> r^{1/n} \cdot m$, except for an explicit finite list of $(r,n)$. }
L. Evain \cite{Ev2} proved this conjecture when the number of points in $\bbP^n$ is of the form $s^n$ (in this case the list is $(4,2)$, $(9,2)$ and $(8,3)$).\\

An affirmative answer to the Nagata conjecture would provide important applications in the theory of linear systems in the projective plane. It would also provide information for 
the study of the singular degrees (works of Bombieri, Skoda, Waldschmidt, Chudnovsky) with a lot of applications in number theory (arithmetic nature of values of Abelian functions of several variables),
in symplectic geometry (symplectic packings in the unit ball) and in algebraic geometry (study of multiple-point Seshadri constants for generic rational surfaces).\\

During the last three/four decades, the use of pluripotential theory in analytic/algebraic geometry has
been very fruitful since it allows a lot of flexibility while keeping track of the analytic features. This
point of view was proven to be very efficient by J.-P. Demailly, Y.-T. Siu and their schools, to name
a few. \\
Our idea in this paper is to develop transcendental techniques, to overcome the intrinsic rigidity of polynomials and to obtain a new approach to this problem of algebraic geometry.
Instead of considering complex polynomials, we work with plurisubharmonic (written psh for simplicity) functions, having logarithmic poles 
at prescribed points. These last functions are much more flexible than the first ones (Lelong). See section $2.1$ for some recalls about pluripotential theory.\\
Thus two points of view are possible. A global one and a local one.\\
We can consider psh functions in all $\bbC^n$, with a logarithmic behavior at infinity.
This class of plurisubharmonic functions contains in particular the logarithm of modulus of polynomials, with prescribed properties in the Nagata conjecture. 
The subclass of plurisubharmonic functions which are maximal (for the complex Monge-Amp\`ere operator+---------) outside prescribed points is of particular interest. 
They inevitably have to satisfy certain conditions of growth at infinity. We look for the minimal growth that they can have. 
In \cite{C-N}, D. Coman and S. Nivoche have obtained some preliminary results about this subject and have established a link between a quantity of the same nature such as the singular 
degrees of M.Waldschmidt \cite{Wal2}, \cite{Wal1}, and a quantity coming directly from pluripotential theory. \\
We can also consider psh functions in a bounded domain in $\bbC^n$, with logarithmic poles at prescribed points and with zero value at the boundary.
We study in particular the subclass of pluricomplex Green functions in this domain with logarithmic poles of weight $1$ at fixed points.
Especially we look at what happens when the poles collide to a single point in the domain. 
The nature of the logarithmic singularity of the limit function is in connection with the algebraic properties of 
the set of fixed points and its singular degree of M.Waldschmidt.\\

In this paper, we establish a link between some well known algebraic quantities of the same nature as the singular 
degree of M.Waldschmidt, and some others quantities coming directly from pluripotential theory. 
We outline new conjectures in term of pluripotential theory, and we prove that they are all equivalent to the Nagata conjecture.

\subsection{Results}

\subsubsection{A first conjecture of pluripotential theory}

We are going to study the convergence of multipole pluricomplex Green functions for a bounded hyperconvex domain in $\bbC^n$, in the case where poles contract to one single point.\\
First we start with the case where the domain is the unit ball centered at the origin. \\
Let $S$ be a finite set of distinct points in $\bbC^n$ and denote $| S |$ its cardinality. Let $R$ be a positive real number sufficiently large such that $S \subset B(O,R)$. 
Denote by $g_R(S,.)$ the pluricomplex Green function in ball $B(O,R)$ 
with logarithmic poles in $S$, of weight one. $g_R(S,z)=g_1(S/R,z/R)$ for any $z \in B(O,R)$. Thus its is natural to study $g_{\infty}$ a negative psh function 
defined in the unit ball $B(O,1)$ by the following upper semi-continuous regularization 
$$g_{\infty}(z)=(\limsup_{t \in \bbC^* \to 0} g_1(tS,z))^*.$$
Clearly, if we replace the set $S$ by $\lambda S$, where $\lambda \in \bbC^*$, then we obtain the same limit function $g_{\infty}$.
$| | .||$ is the Euclidean norm in $\bbC^n$ and $B(z_0,R)$ is the hermitian ball centered at $z_0$ with radius $R$.
$g_{\infty}$ always satisfies in $B(O,1)$ the following inequalities
$$| S |  \cdot \ln || z || \le g_{\infty}(z) \le g_{B(O,1)}(O,z) = \ln || z ||.$$
It tends to $0$ on the boundary of the unit ball and it has an unique logarithmic singularity at the origin.\\
Let us also consider another family $(\tilde{g_1}(tS,.))_{t \in \bbC^*}$ of continuous and psh functions, defined in $B(O,1)$ by 
$\tilde{g_1}(tS,z) = \sup \{g_1(tS,w) : | | w | | = | | z | |\} $ and the following continuous and psh function in $B(O,1)$
$$\tilde{g}_\infty(z) = \sup \{g_ \infty(w) : | | w | | = | | z | | \} .$$
No upper regularization is needed in this definition, since it is a convex increasing (hence continuous) function of $\log|| z | |$ and hence it is itself continuous.
We want to study precisely the nature of the logarithmic singularity of these functions $g_\infty$ and $\tilde{g}_\infty$ at the origin. In particular we want to understand what their Lelong numbers at the origin are.\\

We recall that if $u$ is a psh function, then the {\it classical Lelong number} $\nu(u,z)$ of $u$ at a point $z$ is (P. Lelong, 1969)
the $(2n-2)$-dimensional density of the measure $dd^c u$ at $z$ :
$$\nu(u,z) :=\lim_{r \to 0} \frac{1}{(\pi r^2)^{n-1}} \int_{| w-z| <r} dd^c u \land (dd^c | w-z |^2)^{n-1}.$$
We can also compute this number as follows (V. Avanissian, C. Kiselman):
$$\nu(u,z)=\lim_{y \to -\infty} \frac{\sup_{| w | =1} u(z+we^y)}{y} =
\lim_{y \to -\infty} \frac{1}{y} \int_{| w | =1} u(z+we^y)d\tilde{\lambda}(w),$$
where $d\tilde{\lambda}$ is the normalized surface measure on the unit sphere.\\

There already exist some results about the convergence of multipole Green functions in \cite{M-R-S-T}, \cite{R-T} and \cite{N-T}.\\
In this paper, we establish a direct connection between the nature of the logarithmic singularity of these psh functions $g_\infty$ and $\tilde{g}_\infty$ 
and the algebraic properties of the set $S$. More precisely, by using a Schwarz' Lemma for finite sets (\cite{Mo1}, \cite{Wal1}) and a generalization of a result in \cite{Ni1}, 
we prove in section $2.3$ the following theorem which describes properties of these two psh functions, 
in connection with $\Omega(S)$, the singular degree of $S$ introduced by Waldschmidt (\cite{Wal2}, \cite{Chu1}) and which is an affine invariant in connection with the Nagata Conjecture. \\
For any polynomial $P \in \bbC[z]=\bbC[z_1,\ldots,z_n]$, deg$P$ is its degree and 
$\mbox{ord}(P,p)$ denotes the vanishing order of $P$ at any point $p$. If $l$ is a positive integer we define
$$\Omega(S,l) = \mbox{min}\{\mbox{deg P} : P \in \bbC[z], \mbox{ord}(P,p) \ge l, \forall \, p \in S\}.$$ 
The limit 
$$\displaystyle \Omega(S):=\lim_{l\to +\infty} \Omega(S,l)/l=\inf_{l\ge 1} \Omega(S,l)/l$$ 
exists and is called the singular degree of $S$. \\
On the other hand, according to Theorem $1.1$ (and Example $5.8$) in \cite{R-T}, we already know that the family $(g_1(tS,.))_{t\in\bbC^*}$ converges locally uniformly outside the origin in $B(O,1)$ to $g_{\infty}$.

\begin{Theorem} \label{thm2.1} Let $S$ be a finite set of points in $\bbC^n$. The two psh functions $g_\infty$ and $\tilde{g}_\infty$ satisfy several properties : $\nu(g_{\infty},O) = \Omega(S)$ and
$(dd^c g_{\infty})^n = 0$ in $B(O,1)\setminus\{O\}$. \\
$\Omega(S)^n \le (dd^c g_{\infty})^n(\{O\}) \le \int_{B(O,1)} (dd^c g_{\infty})^n \le | S |$ and we have 
$$g_{\infty}(z) \le \Omega(S) \ln | | z | |, \mbox{ in } B(O,1).$$ 
The family $(\tilde{g}_1(tS,.))_{t\in\bbC^*}$ converges locally uniformly outside the origin in $B(O,1)$ to $\tilde{g}_\infty$ which is equal to $\Omega(S) \ln | | z | |$ 
in $\bar{B}(O,1)$.\\
If $\Omega(S)=| S |^{1/n}$ then  $(dd^c g_{\infty})^n(\{O\})= 
\Omega(S)^n=| S |$ and
$$g_{\infty}(z) = | S |^{1/n} \ln | | z | |,  \, \, \mbox{ in } B(O,1).$$ 
Conversely if $g_{\infty}$ is equal to $\Omega(S) \ln | | . | |$ in $\bar{B}(O,1)$, then 
$\Omega(S)=| S |^{1/n}$.
\end{Theorem}

We will deduce several applications from this theorem, in particular an equivalence between the 
following conjecture of pluripotential theory in $\bbC^n$ and a weak version of the Nagata Conjecture in $\bbP^n$.\\

{\bf Conjecture $({\cal P}_1)$.} {\it In $\bbC^n$, except for a finite number of integer values $r$, for any general set $S=\{p_1, \ldots, p_r\}$ of $r$ points, 
the family of pluricomplex Green functions $(g_{B(O,1)}(tS,.))_{t \in\bbC^*}$ converges locally uniformly outside the origin of $B(O,1)$ to $r^{1/n}g_{B(O,1)}(O,.)$, when $t$ tends to $0$.}\\

{\bf Conjecture $({\cal A}_1)$. Weak Version of the Nagata Conjecture in $\bbP^n$.} {\it In $\bbP^n$, except for a finite number of integer values $r$,  for any general set $S=\{p_1, \ldots, p_r\}$
of $r$ points, $\Omega(S) = r^{1/n}$.}\\

With the homogeneous coordinates $[z_0:z_1:\ldots:z_n]$ in $\bbP^n$, we know that $\bbP^n$ is a complex manifold of dimension $n$, which is covered by $n$ copies 
${\displaystyle U_{i}=\{[z_0:z_1:\ldots:z_n] \in \bbP^n : | z_{i}\neq 0\}}$ of $\bbC^n$. 
If $S=\{p_1, \ldots, p_r\}$ is a set of $r$ distinct points in $\bbP^n$, in the previous definition of $\Omega(S)$, we have just to replace polynomials $P \in \bbC[z]=\bbC[z_1,\ldots,z_n]$ by
homogeneous polynomials in $\bbC[z_0,z_1,\ldots,z_n]$.\\
We can remark that if $n=2$, this Weak version of the Nagata Conjecture is identical to the original one  when the number of points is not a perfect square. More generally in $\bbP^n$, a stronger version 
of the Weak version of the Nagata Conjecture is satisfied for a number of points of the form $s^n$, according to \cite{Ev2}. When the number of points is not of the form $s^n$, this weak version
corresponds to Iarrobino's conjecture in $\bbP^n$.\\

In the previous construction of this function $g_{\infty}$, we can replace the unit ball $B(O,1)$ by any bounded hyperconvex domain $D$ in 
$\bbC^n$ (see the definition in section $2.1$) and the origin by any point $z_o$ in $D$. In this case the function 
$g_{\infty}$ is defined in $D$ by 
$$(\limsup_{t\in \bbC^* \to 0} g_D(z_o+tS,.))^*,$$ 
where $g_D(z_o+tS,.)$ is the pluricomplex Green function in $D$ with logarithmic poles of weight one at any points of the set $z_0 +tS$.
$g_{\infty}$ has similar properties as in Theorem \ref{thm2.1}, where we replace the pluricomplex Green 
function in the unit ball with logarithmic pole at the origin
$\ln | | . | | = g_{B(O,1)}(O,.)$ by the pluricomplex Green function in $D$ with a logarithmic pole at $z_0$, $g_D(z_o,.)$.
First, according to Theorem $1.1$ (and Example $5.8$) in \cite{R-T}, the family $(g_D(z_o+tS,.))_{t\in\bbC^*}$ converges locally uniformly outside $z_o$ in $D$ to $g_{\infty}$. 
In addition, this function satisfies the following theorem.

\begin{Theorem} \label{thm2.1.1} Let $S$ be a finite set of points in $\bbC^n$. Let $D$ be a bounded hyperconvex domain in $\bbC^n$. Fix $z_o$ in $D$.  The psh function $g_\infty$ 
satisfies several properties : \\
$\nu(g_{\infty},z_o) = \Omega(S)$ and
$(dd^c g_{\infty})^n = 0$ in $D \setminus \{z_o\}$.\\
$\Omega(S)^n \le (dd^c g_{\infty})^n(\{z_o\}) \le \int_{D} (dd^c g_{\infty})^n \le | S |$ and we have 
$$g_{\infty}(z) \le \Omega(S) g_D(z_o,z), \mbox{ in } D.$$ 
If $\Omega(S)=| S |^{1/n}$ then  $(dd^c g_{\infty})^n(\{z_o\})= 
\Omega(S)^n=| S |$ and
$$g_{\infty}(z) = | S |^{1/n} g_D(z_o,z),  \, \, \mbox{ in } D.$$ 
Conversely if $g_{\infty}$ is equal to $\Omega(S) g_D(z_o,.)$ in $D$, then 
$\Omega(S)=| S |^{1/n}$.
\end{Theorem}

\subsubsection{New affine invariants and others conjectures of pluripotential theory}

Instead of considering psh functions in bounded domains in $\bbC^n$, here we study a class of entire psh functions in $\bbC^n$, with logarithmic poles in a finite set of points $S$ and with 
a logarithmic growth at infinity. In particular we are interested in the subclass of such functions which are also locally bounded outside of $S$.\\
If $u$ is a psh function in $\bbC^n$, let us denote $\gamma_u$ the following upper limit: 
$$\gamma_u := \limsup_{| z| \to \infty} \frac{u(z)}{\log | | z | |} \in [0,+\infty].$$
If $S = \{p_1,\ldots,p_r\} \subset \bbC^n$ is a finite set of distinct points, we associate to any psh function $u$ in $\bbC^n$, a number $\omega(S,u)$ defined by 
$$\omega(S,u):=\frac{\sum_{j=1}^k \nu(u,p_j)}{\gamma_u},$$
where $\nu(u,z)$ is the Lelong number of the psh function $u$ at a point $z$. For instance, if the psh function $u$ is of the form $\ln | Q |$ where $Q$ is a polynomial, 
then $\nu(u,z)$ is the vanishing order of $Q$ at $z$ and $\gamma_u$ is the degree of $Q$. We associate to $S$ an affine invariant $\omega_{psh}(S)$ defined by 
$$\omega_{psh}(S) := \sup \{ \omega(S,u) : u \in PSH(\bbC^n)\}.$$
We also consider the class of psh functions in $\bbC^n$ which are locally bounded in $\bbC^n \setminus S$. We set 
$$\omega_{psh}^+(S) := \sup \{ \omega(S,u) : u \in PSH(\bbC^n)\cap L^\infty_{loc}(\bbC^n \setminus S)\}.$$
$\omega_{psh}^+(S)$ is also an affine invariant. Then naturally,  a second conjecture of pluripotential theory can be formulated as follow :\\

{\bf Conjecture $({\cal P}_2)$.} {\it In $\bbC^n$, except for a finite number of integer values $r$, for any general set $S=\{p_1, \ldots, p_r\}$ of $r$ points,
$\omega_{psh}(S) = \omega_{psh}^+(S).$}\\

This problem is non-trivial and it is related to the algebraic geometric properties of $S$. There exists another well known invariant of algebraic geometry defined by
$$\omega(S) := \sup \left\{ \frac{\sum_{j=1}^r \mbox{ord}(P,p_j)}{\mbox{deg P}} : P \in \bbC[z]\right\}.$$
$\omega(S)$ is an affine invariant, related to the singular degree of $S$. G.V. Chudnovsky defined in \cite{Chu1} the very singular degree of $S$, $| S |/\omega(S)$. 
Since for any positive integer $l$, $\omega(S) \ge \frac{| S | l}{\Omega(S,l)},$, the following relation between $\omega(S)$ and $\Omega(S)$ is always satisfied : 
$$\omega(S) \ge \frac{| S |}{\Omega(S)}.$$
In section $1.1$, we have stated a general version of the Nagata conjecture (see also Harbourne in \cite{Ha2}) : {\it if $r >9$ then for generic set $S \subset \bbP^2$ 
with $| S | = r$ one has $\sum_{j=1}^r \mbox{ord}(P,p_j) < \sqrt{r} \mbox{deg}P$, for every polynomial $P \in \bbC[z].$}\\
As Nagata has pointed out, this version is equivalent to the uniform one ``the usual Nagata conjecture''. More generally in $\bbP^n$, we state the following conjecture of algebraic geometry.\\

{\bf Conjecture $({\cal A}_2)$.} {\it In $\bbP^n$, except for a finite number of integer values $r$, for any general set $S=\{p_1, \ldots, p_r\}$ of $r$ points,  
$\omega(S) = | S |^{1-\frac{1}{n}}$.}\\

It is well known that conjectures $({\cal A}_1)$ and $({\cal A}_2)$ are equivalent.\\
Finally we can state a last conjecture of pluripotential theory, which can be seen as the dual version of the first one $({\cal P}_1)$. \\

{\bf Conjecture $({\cal P}_3)$.} {\it In $\bbC^n$, except for a finite number of integer values $r$, for any general set $S$ of $r$ points, we have  : 
for any $\epsilon >0$, there exists an entire continuous psh function $v$ in $L^{\infty}_{loc}(\bbC^n \setminus S)$, such that $\nu(v,p) \ge 1$ for any $p \in S$ and 
$\gamma_v \le (1+\epsilon)| S|^{1/n}$.}

\begin{Theorem} \label{thm10.1}
Each conjecture $({\cal P}_1)$, $({\cal P}_2)$ and $({\cal P}_3)$ is equivalent to $({\cal A}_1)$ and $({\cal A}_2)$.
\end{Theorem}

We will prove Theorem \ref{thm10.1} in section $3.3$.\\

Acknowledgments : I would like to thank Professors Andr\'e  Hirschowitz and Joaquim Ro\'e for helpful discussions. 

\section{Conjecture $({\cal P}_1)$ in connection with pluricomplex Green functions in domains in $\bbC^n$}

\subsection{Some recalls of pluripotential theory}

For a bounded domain $D \subset \bbC^n$, the {\it pluricomplex Green function in $D$ with logarithmic poles in a finite subset $S$ of $D$}, generalizes the
one-variable Green functions (for the Laplacian). It is defined by 
$$g_D(S,z) = \sup \{u(z) : u \mbox{ psh on } D, \, u \le 0, \, u(z) \le \ln | | z - p| | + \, O(1) \mbox{ for any point } p \mbox{ in } S\}.$$
(\cite{De1}, \cite{De2}, \cite{Kl3}, \cite{Le1}, \cite{Le2}). \\
If $D$ is hyperconvex (i.e. there exists a continuous psh function  $D \to ]-\infty,0[$) then  
we have an alternative description of the pluricomplex Green functions in terms of the complex Monge-Amp\`ere
operator, namely $g_D(S,.)$ is the unique solution to the following Dirichlet problem:
$$\left\{ \begin{array}{l}
u \mbox{ is plurisubharmonic and negative on } D, \mbox{ continuous on } \bar{D}\setminus S,\\
(dd^c u)^n = 0 \mbox{ on } D \setminus S,\\
u(z) = \ln | | z - p| | +\, O(1) \mbox{ as } z \to p, \, \forall p \in S,\\
u(z) \to 0 \mbox{ as } z \to  \partial D.
\end{array}
\right.$$
In this case, $(dd^c u)^n=\sum_{p \in S} \delta_p$. \\

The complex Monge-Amp\`ere operator is a good candidate to replace Laplacian in one variable. But a important difference between cases of one variable and several variables is that 
this complex Monge-Amp\`ere operator is not linear anymore. \\
Here $d=\partial + \bar{\partial} = \sum_{j=1}^n \frac{\partial}{\partial z_j} + \sum_{j=1}^n \frac{\partial}{\partial \bar{z}_j}$ and $d^c=\frac{i}{2\pi}(\bar{\partial}-\partial)$. $dd^c = \frac{i}{\pi}\partial \bar{\partial}$.
The normalization of the $d^c$ operator is chosen such that we have precisely $(dd^c\ln|| z -p| |)^n=\delta_p$, the Dirac measure at point $p$.  
The $nth$ exterior power of $dd^c = \frac{i}{\pi}\partial \bar{\partial}$ , i.e. $(dd^c)^n = dd^c \wedge \ldots \wedge dd^c$ ($n$ times), defines the complex Monge-Amp\`ere operator in $\bbC^n$. 
If $u \in {\cal C}^2(D)$, then $(dd^c u)^n = \frac{2^n n!}{\pi ^n} \mbox{det}[ \frac{\partial^2 u}{\partial z_j \partial \bar{z}_k} ]dV$, where 
$dV = (\frac{i}{2})^n d z_1 \wedge d \bar{z}_1 \wedge \ldots \wedge d z_n \wedge d \bar{z}_n$ is the usual volume form in $\bbC^n$.\\
It is a positive measure, defined inductively for locally bounded psh functions according to the definition of Bedford-Taylor \cite{B-T1}, \cite{B-T2}, and it can also be extended to psh functions with isolated or compactly supported poles \cite{De4}. 

\subsection{The singular degree of a finite set}

The singular degree of $S$ introduced by Waldschmidt \cite{Wal1} (see also \cite{Chu1}) is an affine invariant in connection with the  Nagata conjecture. If $l$ is a positive integer we define
$$\Omega(S,l) = \mbox{min}\{\mbox{deg P} : P \in \bbC[z], \mbox{ord}(P,p_j) \ge l, 1\le j \le r\}.$$ 
$\Omega(S,1)$ is sometimes called the degree of $S$. We clearly have $\Omega(S,l_1+l_2) \le \Omega(S,l_1) + \Omega(S,l_2)$, and in particular $\Omega(S,l) \le \Omega(S,1)l$. The limit 
$$\Omega(S):=\lim_{l\to +\infty} \frac{\Omega(S,l)}{l}$$
exists and is called the singular degree of $S$. We have for all $l\ge 1$
$$\frac{\Omega(S,1)}{n} \le \Omega(S) = \inf_l \frac{\Omega(S,l)}{l} \le \frac{\Omega(S,l)}{l} \le \Omega(S,1).$$
The second and the third inequalities are trivial, while the proof of the first one uses complex analysis. By using Hormander-Bombieri-Skoda theorem, 
M. Waldschmidt proved more generally that for any positive integer $l_1$ and $l_2$
$$\frac{\Omega(S,l_1)}{l_1+n-1}\le \Omega(S) \le \frac{\Omega(S,l_2)}{l_2}.$$
Upper bound for the numbers $\Omega(S,l)$ is well known. It is a result of Waldschmidt (\cite{Wal1}, Lemma 1.3.13): 
$$\Omega(S,l) \le (l+n-1) | S |^{1/n} -(n-1).$$
And a consequence is $\Omega(S) \le | S |^{1/n}.$ An important and difficult problem is to find a lower bound for $\Omega(S)$.
The Nagata conjecture can be stated again in term of the invariants $\Omega(S,l) $ : \\
{\it In $\bbC^2$, if $r >9$, then $\Omega(S,l) >l\sqrt{r}$ $\forall \, l \ge 1$ holds for a set $S$ of $r$ points in general position.}\\
This statement doesn't hold for $r \le 9$.

\subsection{Proof of Theorem \ref{thm2.1}}

Let $S$ be a finite set of distinct points in $\bbC^n$. We use the same notation as in section $1.2.1$.
For any $t \in \bbC^*$ sufficiently small, the function $g_1(tS,.)$ always satisfies the following inequalities
$$\sum_{p\in S} g_1(tp,z) \le g_1(tS,z) \le \inf_{p \in S} g_1(tp,z),$$
where $g_1(tp,.)$ is the pluricomplex Green function with an unique logarithmic pole at $tp$ with weight one. Consequently, the psh function $g_{\infty}$ satisfies in $B(O,1)$ 
$$| S | . \ln || z || \le g_{\infty}(z) \le g_1(O,z) = \ln || z ||.$$
In addition, this function tends to $0$ on the boundary of the unit ball and it has an unique singularity at the origin.\\
The functions $g_1(tS,.)$ are continuous in $\bar{B}(O,1)$ valued in $[-\infty,0]$ and psh in $B(O,1)$ with a finite set of singularities $tS$. 
No upper regularization is needed in the definitions of $\tilde{g_1}(tS,.)$ and $\tilde{g}_\infty$, because they are convex increasing (hence continuous) functions of $\log|| z | |$ and hence they are themself continuous.\\

The following very simple example, helps us to understand what different situations can occur in Theorem \ref{thm2.1}.

\begin{Example}
Let $S=\{(1/2,0),(-1/2,0)\}$ in $\bbC^2$. We know that $\Omega(S)=1$. \\
To simplify computation, we use the sup norm $| z | = \mbox{max}\{| z_1 |, | z_2 |\}$ instead of the Euclidean norm and we are in the unit polydisc instead of the unit ball. For any $t \in \bbC^*$ sufficiently small, 
$$g_1(tS,z)=\mbox{max}\{\ln | \frac{(z_1 - t/2)(z_1 + t/2)}{(1-\bar{t}z_1/2)(1+\bar{t}z_1/2)} |, \ln | z_2 |\}$$ 
is the pluricomplex Green function in the unit polydisc $P(O,1)$ with logarithmic poles in $tS$ with weight $1$. $(g_1(tS,.))_{t \in \bbC^*}$ converges uniformly on any compact set of the form $\overline{P(O,1)} \setminus P(O,\varrho)$ (where $0< \varrho <1$) to $g_{\infty}$ which is explicitely equal to $g_\infty(z)=\mbox{max}\{2\ln | z_1 |, \ln | z_2 |\}$ 
in $P(O,1)$. It is also easy to verify that $\tilde{g}_\infty(z)=\mbox{max}\{\ln | z_1 |, \ln | z_2 |\}=\ln | z |$ in $P(O,1)$. 
\end{Example}

\noindent  {\sc Proof of Theorem \ref{thm2.1}.} \\
In the first three sections and section $7.$, we prove item $(ii)$. In sections $4.$ and $5.$, we prove item $(i)$. 
In section $6.$, we prove item $(iii)$ and in section $8.$, we prove item $(iv)$.\\

\indent {\bf 1)} Let $Q$ be a holomorphic polynomial and $\varrho$ be a positive real number. It is well known that 
$$\frac{1}{deg(Q)} \ln \left(\frac{| Q(z) |}{| | Q | |_\varrho}\right) \le \ln^+ 
\left(\frac{| | z | |}{\varrho}\right), \, \mbox{ in } \bbC^n.$$
$|| f | |_\varrho = \sup_{z \in B(O,\varrho)} | f(z) |$ is the sup-norm in $B(O,\varrho)$. The function on the right hand side, is called the pluricomplex Green function with 
logarithmic pole at infinity for the compact set $\bar{B}(O,\varrho)$ (Zahariuta, Siciak). $\ln^+$ is the function $\mbox{max}\{\ln, 0\}$. 
Consequently, we obtain for any positive real number $R$ such that $0<\varrho \le R$,
\begin{equation} \label{eq1}
\ln || Q | |_\varrho - \ln | | Q | |_R \, \ge \, deg(Q).\ln \left(\frac{\varrho}{R}\right).
\end{equation}
On the other hand, according to a Schwarz' Lemma (Moreau \cite{Mo1}, Waldschmidt \cite{Wal1} p.146): for any $\epsilon >0$, there exists a real number $r_2(S,\epsilon)$ such that for any 
$r_2 < \varrho < R$ with $2e^n.\varrho < R$ and for any polynomial $Q$ such that $\mbox{ord}(Q,p) \ge l$ for all $p \in S$,
\begin{equation} \label{eq2}
\ln || Q | |_\varrho - \ln | | Q | |_R \, \le \, (\Omega(S,l) - l\epsilon).\ln \left(\frac{2 e^n \varrho}{R}\right) \le -
l(\Omega(S)-\epsilon)\ln \left(\frac{R}{2 e^n \varrho}\right).
\end{equation}
We can suppose in addition that $r_2(S,\epsilon) \ge \sup_{p \in S} | | p| | := | |  S| |$.\\

{\bf 2)} We can generalize Theorem $1.1$ of \cite{Ni1} (for a pluricomplex Green function with one logarithmic pole) to the case of a pluricomplex Green function with a finite number
of logarithmic poles: 
$$g_{R}(S,.) = \sup_{l \ge 1} H_{S,R,l}=\lim_{l\to \infty}  H_{S,R,l},$$ 
where each psh function $H_{S,R,l}$ ($l$ is a positive integer) is defined by 
$$H_{S,R,l} = \sup \{\frac{1}{l} \ln | f |  : f \in {\cal O}(B(O,R)), \, | | f | |_{R} \le 1,   \mbox{ ord}(f,s) \ge l  \mbox{ for any } s \in S\}.$$
Here ${\cal O}(D)$ is the Frechet space of holomorphic functions in an open set $D$.\\
We can replace in the previous definition of $H_{S,R,l}$, holomorphic functions by polynomials satisfying the same properties. 
Indeed, if $l$ and $d$ are two positive fixed integers, then we have a continuous linear map $\varphi$ from $\bbC_d [z]$ to $\bbC^{m_l}$ ($m_l = | S| \binom{l-1+n}{n}$) as follow
$$\varphi : P \mapsto (P^{(\alpha)}(s), \forall | \alpha | \le l-1, \, \forall s \in S).$$
$\bbC_d [z]$ is the linear space of polynomials of degree less or equal to $d$ of dimension $\binom{d+n}{n}$. If $d$ is sufficiently large, 
then this map $\varphi$ is surjective (we can prove it  by using Cartan-Serre's Theorem). $\mbox{Ker} \varphi$ is the kernel of $\varphi$ and it has 
dim$(\mbox{Ker}  \varphi)=\binom{d+n}{n} - | S| \binom{l-1+n}{n}$. 
We choose $E$ a subspace of $\bbC_d [z]$, such that $\mbox{Ker}  \varphi \oplus E = \bbC_d [z]$. Then there exists a map $\psi$ from $\bbC^{m_l}$ to $E$, 
which is the inverse of the restriction of $\varphi$ to $E$. $\psi$ is a bijective continuous linear map from $\bbC^{m_l}$ to $E$ with norm $|| \psi ||$ (which depends on $l$ and $d$).\\
According to the fact that the family $(H_{S,R',l})_{R'>R}$ converges uniformly on $B(O,R)$ to  $H_{S,R,l}$ when $R'$ goes to $R$ and because $B(O,R)$ is a Runge domain in $\bbC^n$; for any $f$ with $| | f | |_{R} = 1$ which appears in the definition of $H_{S,R,l}$, we can find a sequence $(P_j)_j$ of polynomials such that $| | P_j | |_{R} \le 1$ and 
$| | f-P_j | |_{R} \le 1/j$.\\
There exists $\delta >0$ such that $P(s,\delta) \subset B(O,R)$, for all $s \in S$. $P(s,\delta)$ is the polydisc centered in $s$ with multiradius $(\delta,\ldots,\delta)$. According to Cauchy's inequalities applied in each polydisc $P(s,\delta)$ and since 
$| | f-P_j | |_{P(s,\delta)} \le 1/j $, we obtain that $| (f-P_j)^{(\alpha)}(s) | = | (P_j)^{(\alpha)}(s) | \le \alpha!/(jr^{| \alpha |})$, 
for all $| \alpha | \le l-1$ and for all $s \in S$.\\
Denote by $Q_j$ the polynomial in $\bbC_d [z]$, equal to $\psi((P_j)^{(\alpha)}(s), \forall | \alpha | \le l-1, \, \forall s \in S)$. There exists a positive constant $c_{l,d}$ such that 
$| | Q_j | |_{R} \le c_{l,d}/j$. Then if we replace $P_j$ by $R_j:=(P_j - Q_j)/| | P_j -Q_j | |_{R}$ we obtain that $| | f-R_j | |_{R} \le 1+O(c_{l,d}/j)$, with 
ord($R_j, p) \ge l$, for all $p \in S$.\\

{\bf 3)} Since for any positive integer $l $, there exists a polynomial $Q$ such that $\mbox{ord}(Q,p) \ge l$ for any $p \in S$, $deg(Q)=\Omega(S,l)$ and $| | Q | |_{R} =1$, according to inequalities $(\ref{eq1})$ and $(\ref{eq2})$, we obtain 
$$ \frac{\Omega(S,l)}{l} \ln \left(\frac{\varrho}{R}\right) \le \sup_{z \in B(O,\varrho)} H_{S,R,l}(z) \le (\frac{\Omega(S,l)}{l}-\epsilon) \ln \left(\frac{2 \varrho e^n}{R}\right).$$
Since $\displaystyle H_{S,1/| t|,l}( z) = H_{tS,1,l}(tz)$ when $| z | \le 1/| t|$ and $t \in \bbC^*$,   we have in particular 
$$ \frac{\Omega(S,l)}{l} \ln (| t | \varrho) \le \sup_{z \in B(O,| t |\varrho)} H_{tS,1,l}(z) \le (\frac{\Omega(S,l)}{l}-\epsilon) \ln (2 \varrho | t | e^n),$$
or which is equivalent: for any $t \in \bbC^*$ and $\varrho' >0$ such that $\displaystyle  | t | r_2(S,\epsilon) \le \varrho' \le 1$
$$ \frac{\Omega(S,l)}{l} \ln \varrho' \le \sup_{z \in B(O,\varrho')} H_{tS,1,l}(z) \le (\frac{\Omega(S,l)}{l}-\epsilon) \ln \left(2 \varrho' e^n \right).$$
We can suppose in addition that $t$ is sufficiently small such that $tS \subset B(O,\varrho')$. We know that the sequence $(H_{tS,1,l})_l$ converges uniformly 
on $\bar{B}(O,1)\setminus B(O,\varrho')$ to $g_{1}(tS,.)$ when $l$ tends to infinity (\cite{Ni1}). From the previous estimates, we deduce 
$$ \Omega(S) \ln \varrho' \le \sup_{z \in B(O,\varrho')} g_1(tS,z)  \le (\Omega(S) -\epsilon) \ln \left(2 \varrho' e^n \right),$$
for any $\varrho' >0$ and $t \in \bbC^*$ satisfying $| t | r_2(S,\epsilon) \le \varrho' \le 1$ and $tS \subset B(O,\varrho')$.\\
If $z \ne O$, $| t | r_2(S,\epsilon) \le | | z | | \le 1$ and $tS \subset B(O,| | z | |)$, 
\begin{equation} \label{eq3}
\Omega(S) \ln | | z | | \le \tilde{g}_1(tS,z)  \le (\Omega(S) -\epsilon) \ln \left(2 | | z | | e^n \right).
\end{equation}
In particular if $| t | r_2(S,\epsilon) = | | z | |$ (we have $r_2(S,\epsilon) \ge | |  S| |$), we have 
$$\tilde{g}_1(tS,z)  \le (\Omega(S) -\epsilon) 
\ln ( 2 e^n | t | r_2(S,\epsilon))= \left( (\Omega(S) -\epsilon)\frac{\ln (2e^n)}{\ln(| t | r_2)} + (\Omega(S) -\epsilon) \right) \ln  | | z | |.$$
For $t$ sufficiently small, $(\Omega(S) -\epsilon)\frac{\ln (2e^n)}{\ln(| t | r_2)} + (\Omega(S) -\epsilon) >0$, the psh function \\$\left( (\Omega(S) -\epsilon)\frac{\ln (2e^n)}{\ln(| t | r_2)} + (\Omega(S) -\epsilon) \right) \ln | | . | |$ is maximal in $B(O,1)  \setminus \{O\}$ and equal to $0$ on $\partial B(O,1)$. Thus we obtain in $B(O,1)  \setminus B(O,| t | r_2)$,
$$\Omega(S) \ln | | z | | \le \tilde{g}_1(tS,z)  \le \left( (\Omega(S) -\epsilon)\frac{\ln (2e^n)}{\ln(| t | r_2)} + (\Omega(S) -\epsilon) \right) \ln | | z | |.$$
Consequently, for any $0<\varrho <1$, $(\tilde{g}_1(tS,.))_t$ converges uniformly in $\bar{B}(O,1)  \setminus B(O,\varrho)$ to 
$\Omega(S) \ln | | . | |$ when $t$ converges to $0$.\\

{\bf 4)} Let $0< \varrho <1$ be fixed. Let $\epsilon_1$ be a positive real number such that $\epsilon_1 r_2(S,\epsilon) \le \varrho \le 1$. Since $\sup_{| t | \le \epsilon_1} \sup_{z \in B(O,\varrho)} g_1(tS,z) = \sup_{z \in B(O,\varrho)} \sup_{| t | \le \epsilon_1} g_1(tS,z)$, we obtain according to the previous estimates $(\ref{eq3})$
\begin{equation} \label{eq4}
\Omega(S) \ln \varrho \le \sup_{z \in B(O,\varrho)} (\sup_{| t | \le \epsilon_1}  g_1(tS,z))^*  \le (\Omega(S) -\epsilon) \ln \left(2\varrho e^n \right).
\end{equation}
Since $\displaystyle g_{\infty} = \lim_{\epsilon_1 \to 0} (\sup_{| t | \le \epsilon_1}  g_1(tS,.))^*$, where this limit decreases when $\epsilon_1$ decreases, and according to the fact that the family $(g_1(tS,.))_{t \in \bbC^*}$ converges locally uniformly outside the origin in $B(O,1)$ to $g_{\infty}$ (Theorem $1.1$ of \cite{R-T}), we can deduce that
$$\Omega(S)  \ln \varrho \le \sup_{z \in B(O,\varrho)} g_{\infty}(z)  \le (\Omega(S) -\epsilon) \ln \left(2 \varrho e^n \right).$$
Consequently, $\displaystyle \lim_{\varrho \to 0} \frac{\sup_{z \in B(O,\varrho)} g_{\infty}(z)}{\ln \varrho} := \nu(g_{\infty},O) =\Omega(S)$. In particular, according to the maximality of the psh function $\ln | | z | |$ in $B(O,1)\setminus \{O\}$, we obtain in $B(O,1)$ that
$$g_{\infty}(z) \le \Omega(S) \ln | | z | |.$$
And $(dd^c g_{\infty})^n(\{O\}) \ge \Omega(S)^n$.\\

{\bf 5)} By the comparison principle for the complex Monge-Amp\`ere operator \cite{B-T2}, we have for any $\epsilon >0$
$$\int_{\{u<v\}} (dd^c v)^n \le \int_{\{u<v\}} (dd^c u)^n,$$
where $u$ and $v$ are two psh functions in $B(O,1)$ defined by $u=(1+\epsilon)g_1(tS,.)$ and $v=(\sup_{| t | \le \epsilon_1}  g_1(tS,.))^*$. $u$ has isolated logarithmic singularities and $v$ is bounded. Since $u$ and $v$ tend to $0$ on the boundary 
of the unit ball, we deduce that 
$$\int_{B(O,1)} (dd^c (\sup_{| t | \le \epsilon_1}  g_1(tS,.))^*)^n \le \int_{B(O,1)} (dd^c (1+\epsilon)g_1(tS,.))^n=(1+\epsilon)^n | S |.$$
According to the monotonicity of the family $((\sup_{| t | \le \epsilon_1}  g_1(tS,.))^*)_{\epsilon_1}$ when $\epsilon_1$ decreases, and by making tend $\epsilon$ to $0$, we deduce that 
$$\int_{B(O,1)} (dd^c g_{\infty})^n \le  | S |.$$
Since the family $(g_1(tS,.))_{t \in \bbC^*}$ converges locally uniformly outside the origin in $B(O,1)$ to $g_{\infty}$, then the complex Monge-Amp\`ere measure
$(dd^c g_1(tS,.))^n$ converges to $(dd^c g_{\infty})^n$ in $B(O,1)\setminus \{O\}$. Since $(dd^c g_1(tS,.))^n=\sum_{p\in S} \delta_{tp}$, we have $(dd^c g_{\infty})^n=0$ in $B(O,1)\setminus \{O\}$. \\

{\bf 6)} If $\Omega(S)=| S |^{1/n}$ then all inequalities in item $(i)$ are equalities. $\nu(g_{\infty},O) = \Omega(S)$, $(dd^c g_{\infty})^n = 0$ in $B(O,1)\setminus\{O\}$ and $(dd^c g_{\infty})^n(\{O\})= \Omega(S)^n=| S |$.\\ 
We conclude this case by proving that $g_{\infty}(z) = | S |^{1/n} \ln | | z | |$ in $B(O,1)$.\\
Assume that there exists a point $a \in B(O,1)\setminus\{O\}$ such that $g_{\infty}(a) < | S |^{1/n} \ln | | a | |.$ Let $\gamma$ be a ${\cal C}^{\infty}(\overline{B(O,1)})$ and 
strictly psh function such that $\gamma \le -1$ in $\overline{B(O,1)}$ (classical argument, used in the proof of Theorem $4.3$ \cite{De2} for instance. We can choose $\gamma(z)=|| z ||^2 -2$). Let $\epsilon > 0$ be sufficiently small such that the following psh function $v$ defined by
$$v(z)=\mbox{max}\{\Omega(S) \ln | | z | | + \epsilon \gamma(z),g_{\infty}(z)\}$$
satisfies\\
{\it (i)} $v(z)= \Omega(S) \ln | | z | | + O(1)$ when $z \to O$,\\
{\it (ii)} $v(z)=g_{\infty}(z)$ in a neighborhood of $\partial B$,\\
{\it (iii)} $v(z)= \Omega(S) \ln | | z | | + \epsilon \gamma(z)$ in a neighborhood of $a$.\\
According to $(i)$ and $(ii)$, we deduce that ${\bf 1}_{\{O\}}(dd^cv)^n = \Omega(S)^n \delta_O= | S | \delta_O$ and 
$\int_{B(O,1)} (dd^c v)^n $ $ = \int_{B(O,1)} (dd^c g_{\infty})^n = | S |$ 
(according to the comparison principle for MA operator). Hence $(dd^c v)^n=0$ in $B(O,1)\setminus\{O\}$. This is in contradiction with the third above property. Consequently $g_{\infty} \equiv | S |^{1/n} \ln | | . | |$ in $B(O,1)$.\\

{\bf 7)} 
According to the third section of this proof, $g_\infty (z) \le \Omega(S) \ln | | z | |$ in $\bar{B}(O,1) $. Then $\tilde{g}_\infty$ satisfies the same estimate in $\bar{B}(O,1)$ :
$$\tilde{g}_\infty (z) \le \Omega(S) \ln | | z | |.$$
Let us prove that $\tilde{g}_\infty (z) = \Omega(S) \ln | | z | |$ in $\bar{B}(O,1) $. Suppose that there exists $O \neq z \in B(O,1)$, such that $-\infty < \tilde{g}_\infty (z) =
c \ln | | z | | < \Omega(S) \ln | | z | |$. Then $\tilde{g}_\infty (w)$ satisfies the same estimates for any $w$ such that $|| w || = || z ||$ and we deduce that
$$g_\infty(w) \le c \ln || w || < \Omega(S) \ln || w ||, \, \forall w \in \partial B(O,|| z ||).$$
By definition we have $\displaystyle \limsup_{t\to 0} g_1(tS,w) \le (\limsup_{t\to 0} g_1(tS,w))^* = g_\infty(w)$. According to Hartog's Lemma, we obtain for the positive constant $d= (\Omega(S) - c) \ln || z ||$, that there exists $\epsilon_0$ such that for any $| t | \le \epsilon_0$
$$g_1(tS,w) \le c \ln || w || + d/2 < \Omega(S) \ln || w ||,\, \forall w \in \partial B(O,|| z ||).$$
In particular, $\displaystyle \sup_{w \in \partial B(O,|| z ||)} g_1(tS,w) \le c \ln || z || + d/2 < \Omega(S) \ln || z ||$. This contradicts what it has been proved in section $3.$ and consequently, $\tilde{g}_\infty \equiv \Omega(S) \ln | | . | |$ in $\bar{B}(O,1) $. \\

{\bf 8)} Let us prove the last item {\it (iv)}. We suppose that $(g_1(tS,.))_{t \in \bbC^*}$ converges locally uniformly outside the origin in $B(O,1)$ to $g_{\infty}=\Omega(S) \ln || . ||$ in $\bar{B}(O,1)$. Then for any $\epsilon >0$ and for any $0<\varrho <1$, there exists $\epsilon_0>0$ such that for any $| t | \le \epsilon_0$,
we have $tS \subset B(O,\varrho)$ and
$$(1+\epsilon) \Omega(S) \ln || z || \le g_1(tS,z) \le (1-\epsilon) \Omega(S) \ln || z ||, \, \forall z \in \bar{B}(O,1)  \setminus B(O,\varrho).$$
Let us construct a continuous and psh entire function $v$ such that 
$$v(z) = \left\{ \begin{array}{ll}
g_1(tS,z), & z \in B(O,\varrho),\\
\mbox{max}\{g_1(tS,z),(1+\epsilon) \Omega(S) \ln || z ||\},  & z \in B(O,1) \setminus B(O,\varrho),\\
(1+\epsilon) \Omega(S) \ln || z || , & z \in \bbC^n \setminus B(O,1).          
\end{array}\right.$$
We obtain that $| S | = \int_{B(O,\varrho)} (dd^c v)^n \le \int_{\bbC^n} (dd^c v)^n = (1+\epsilon)^n \Omega(S)^n$. We already know that 
$\Omega(S) \le | S |^{1/n}$. Consequently, $\Omega(S) = | S |^{1/n}$. This achieves the proof of Theorem \ref{thm2.1}.\cqfd\\


\noindent  {\sc Sketch of the proof of Theorem \ref{thm2.1.1}.} \\
$D$ is a bounded hyperconvex domain in $\bbC^n$ and $z_o$ is any point in $D$. There exists two positive real numbers $R_1$ and $R_2$ such that $B(z_0,R_1) \subset D \subset B(z_0,R_2)$. Then 
$$g_{B(z_0,R_2)}(z_0+tS,z) \le g_{D}(z_0+tS,z) \le g_{B(z_0,R_1)}(z_0+tS,z) \,\,\,\mbox{ in } B(z_0,R_1)$$ 
and 
$$g_{B(z_0,R_i)}(z_0+tS,z) = g_{B(O,1)}(\frac{tS}{R_i},\frac{z-z_0}{R_i}) \,\,\,\mbox{ in } B(z_0,R_i), \,\,\,\mbox{ for } i=1,\, 2.$$
Thus we have
$$g_{B(O,1)}(\frac{tS}{R_2},\frac{z-z_0}{R_2}) \le g_{D}(z_0+tS,z) \le g_{B(O,1)}(\frac{tS}{R_1},\frac{z-z_0}{R_1}),  \mbox{ in }
B(z_0,R_1).$$
$g_{\infty}$ is defined by 
$$g_{\infty} = (\limsup_{t\to 0} g_D(z_0+tS,.))^*, \,\,\,\mbox{ in } D.$$
According to Theorem $1.1$ (and Example $5.8$) in \cite{R-T}, the family $(g_D(z_0+tS,.))_{t\in\bbC^*}$ converges locally uniformly outside the point $z_o$ in $D$ to $g_{\infty}$.\\
Let us denote by $g_{\infty,i}$ the following psh function
$$g_{\infty,i} = (\limsup_{t\to 0} g_{B(O,1)}(tS,.))^*, \,\,\,\mbox{ in } B(O,1),\,\,\,\mbox{ for } i=1,\, 2.$$
Since $\nu(g_{\infty,i},O)= \Omega(S)$ for $i=1,\, 2$, $\nu(g_{\infty},z_0) =\Omega(S)$,
$g_{\infty} \le \Omega(S) g_{D}(z_0,.) \mbox{ in } D$ and $(dd^c g_{\infty})^n(\{z_o\}) \ge \Omega(S)^n$.\\
As in section $5.$ of the previous proof, we obtain that $(dd^c g_{\infty})^n = 0$ in $D \setminus \{z_o\}$ and that $\Omega(S)^n \le \int_{D} (dd^c g_{\infty})^n \le | S |$. Item $(i)$ of Theorem \ref{thm2.1.1} is proved.\\
To prove item $(ii)$, we proceed as in section $6.$ of the previous proof. $\gamma$ is replaced by a ${\cal C}^{\infty}(\overline{D})$ and strictly psh function such that $\gamma \le -1$ in $\overline{D}$ (classical argument, used in the proof of Theorem $4.3$ \cite{De2} for a bounded hyperconvex domain in $\bbC^n$. We can choose $\gamma(z)=(|| z -z_o||/R_2)^2 -2$).\\
To prove the last item {\it (iii)}, we proceed as in section $8.$ of the previous proof. The family $(g_D(z_0+tS,.))_{t\in\bbC^*}$ converges locally uniformly in $D \setminus \{z_0\}$ to $g_{\infty}=\Omega(S) g_D(z_o,.)$ when $t$ goes to $0$.\\
Then for any $\epsilon >0$ and for any $\varrho > 0$, there exists $\epsilon_0>0$ such that for any $| t | \le \epsilon_0$,
we have $tS \subset B(O,\varrho)$ and
$$(1+\epsilon) \Omega(S) g_D(z_o,z) \le g_D(z_o+tS,z) \le (1-\epsilon) \Omega(S) g_D(z_o,z), \, \forall z \in \bar{D}  \setminus B(z_o,\varrho).$$
Let us construct a continuous and psh function $v$ in $D$ such that 
$$v(z) = \left\{ \begin{array}{ll}
g_D(z_o+tS,z), & z \in B(z_o,\varrho),\\
\mbox{max}\{g_D(z_o+tS,z),(1+\epsilon) \Omega(S) g_D(z_o,z)\},  & z \in D \setminus B(z_o,\varrho).        
\end{array}\right.$$
We obtain that $| S | = \int_{B(z_o,\varrho)} (dd^c v)^n \le \int_{D} (dd^c v)^n \le (1+\epsilon)^n \Omega(S)^n$. We already know that 
$\Omega(S) \le | S |^{1/n}$. Consequently, $\Omega(S) = | S |^{1/n}$. This achieves the proof of Theorem \ref{thm2.1.1}.\cqfd

\subsection{A Schwarz' Lemma}


With the first item of the previous theorem, we obtain a slightly different version of Schwarz' Lemma from the previous ones \cite{Mo1}. 

\begin{Corollary} \label{cor2.1}
For any positive real numbers $\varrho <1$ and $\epsilon$, there exits a positive real number $R_{\epsilon, \varrho}$ such that  : 
if $R$ and $r$ satisfy $R \ge R_{\epsilon, \varrho}$ and $\varrho R \le r \le R$ we have :
$$\ln | | f | |_r \le \ln | | f | |_R - l(\Omega(S) - \epsilon) \ln (R/r),$$
for any integer $l \ge 1$ and for any entire function $f$ with ord$(f,p) \ge l$ for all $p \in S$.
\end{Corollary}

{\sc proof.} According to Theorem \ref{thm2.1}, $g_\infty(z) \le \Omega(S) \ln | | z ||$ in $B(O,1)$. If we apply Hartogs lemma to the compact set $\partial B(O,\varrho)$, we obtain that: there exists a positive real number $R_{\epsilon, \varrho}$ such that  for any $R \ge R_{\epsilon, \varrho}$, we have $S/R \subset B(O,\varrho)$ and
$$g_1(S/R,z) \le (\Omega(S)-\epsilon) \ln | | z | | , \, \, \mbox{ on } \partial B(O,\varrho).$$
$(\Omega(S)-\epsilon) \ln | | . | |$ is maximal in $B(O,1)\setminus \{O\}$ and these two functions are equal to $0$ on the boundary of $B(O,1)$. Consequently, the same inequality is satisfied in $B(O,1) \setminus B(O,\varrho)$ and
$$g_R(S,z) \le (\Omega(S)-\epsilon) \ln | | z/R | | , \, \, \mbox{ on } B(O,R) \setminus B(O,\varrho R).$$
Finally, if $l$ is a positive integer and $f$ is an entire function with ord$(f,p) \ge l$ for all $p \in S$ ($\ln (| f | /| | f |  |_R) \le l.g_R(S,.)$ in $B(O,R)$), then for any $\varrho R \le r \le R$ we obtain 
$$\ln | | f | |_r \le \ln | | f | |_R - l(\Omega(S) - \epsilon) \ln (R/r).$$\cqfd 

\subsection{Pluricomplex Green functions with a finite set of logarithmic poles with different homogeneous weights}

Let $D$ be a bounded hyperconvex domain in $\bbC^n$. Let $S=\{p_1, \ldots, p_r\}$ be a finite set of $r$ distinct points in $D$ and ${\cal M}=\{m_1,\ldots, m_r\}$ be a set of $r$ positive integers. 
Denote by $g_D((S,{\cal M}),.)$, the pluricomplex Green function in the domain $D$, with logarithmic poles at any points $p_j$ of $S$ and with homogeneous weight $m_j$ respectively.\\
According to \cite{Ev2}, to any point $p_j \in S$ we can associate a generic set $S_j$ of $m_j^n$ distinct points in $\bbC^n$ such that $\Omega(S_j)=m_j$. \\
Denote by $\Sigma(t_1,\ldots,t_r)$ the set $\{p_j+t_j S_j : 1 \le j \le r\}$ of $\sum_{j=1}^r m_j^n$ distinct points. For sufficiently small $(t_1,\ldots,t_r) \in (\bbC^*)^r$ the set $\Sigma(t_1,\ldots,t_r)$ is contained in $D$. Denote by $g_D(\Sigma(t_1,\ldots,t_r),.)$ the pluricomplex Green function in $D$, with logarithmic poles at any points of $\Sigma(t_1,\ldots,t_r)$ and with weight $1$. These functions satisfy in $D$
$$\sum_{j=1}^{r} g_D(p_j+t_j S_j,z) \le g_D(\Sigma(t_1,\ldots,t_r),z) \le \inf_{1\le j \le r} g_D(p_j+t_j S_j,z).$$
According to Theorem \ref{thm2.1}, for any $1\le j \le r$, the family $(g_D(p_j+t_j S_j,.))_{t_j\in\bbC^*}$ converges locally uniformly outside the point $p_j$ in $D$ to $m_j g_D(p_j,.)$. \\
By using the same technics as in the proof of Theorem \ref{thm2.1}, we deduce the following theorem.

\begin{Theorem} \label{thm3.1}
The family of pluricomplex Green functions $(g_D(\Sigma(t_1,\ldots,t_r),.))_{(t_1,\ldots,t_r) \in (\bbC^*)^r}$ converges locally uniformly outside the set $S$ in $D$ to $g_D((S,{\cal M}),.)$, when $(t_1,\ldots,t_r)$ tends to $O$.
\end{Theorem}

\section{New affine invariants and others conjectures of pluripotential theory}


\subsection{Affine invariants associated to finite sets of points in $\bbC^n$}

Fix a finite set $S = \{p_1,\ldots,p_r\} \subset \bbC^n$ of $r$ distinct points, $| S | = r$ is its cardinality.\\
In section 1.2.2 we have defined two affine invariants $\omega_{psh}(S)$ and $\omega_{psh}^+(S)$. 
$$\omega_{psh}(S) := \sup \{ \omega(S,u) : u \in PSH(\bbC^n)\}$$
and 
$$\omega_{psh}^+(S) := \sup \{ \omega(S,u) : u \in PSH(\bbC^n)\cap L^\infty_{loc}(\bbC^n \setminus S)\},$$
where for any psh function $u$ in $\bbC^n$, $\displaystyle \gamma_u := \limsup_{| z| \to \infty} \frac{u(z)}{\log | z|}$ 
and $\displaystyle \omega(S,u):=\frac{\sum_{j=1}^k \nu(u,p_j)}{\gamma_u}$. $\nu(u,z)$ is the Lelong number of the psh function $u$ at a point $z$. \\
Let us remark that for any psh function $u$ in $\bbC^n$, we have always $\gamma_u \in [0,+\infty]$. Indeed, if $u$ is psh in $\bbC^n$ and $u(z)\le o(\log | z |)$ as $| z | \to \infty$, then $u$ is constant 
(it is a generalization of Liouville's Theorem, see for example \cite{Ho1}) and $\gamma_u = 0$. 
Conversely, if $u(z)=| z |$ then $\gamma_u=+\infty$.\\
$\omega_{psh}(S)$ and $\omega_{psh}^+(S)$ are related to the algebraic geometric properties of $S$ and in particular to the affine invariant $\omega(S)$ of $S$, that they generalize. \\
Here are some simple properties of these invariants.


\begin{Lemma}\label{lm1.1} 
{\it (i)} $\omega(S) \le \omega_{psh}(S)$ and $\omega_{psh}^+(S) \le \omega_{psh}(S).$ \\
{\it (ii)} The sets $\{ \omega(S,u) : u \in PSH(\bbC^n)\}$ and $\{ \omega(S,u) : u \in PSH(\bbC^n)\cap L^\infty_{loc}(\bbC^n \setminus S)\}$ are connected.
$$(0,\omega_{psh}(S))\subseteq \{\omega(S,u) : u \in PSH(\bbC^n) \} \subseteq (0,\omega_{psh}(S)],$$ 
$$(0,\omega_{psh}^+(S))\subseteq \{\omega(S,u) : u \in PSH(\bbC^n)\cap L^\infty_{loc}(\bbC^n \setminus S) \} 
\subseteq (0,\omega_{psh}^+(S)].$$
\end{Lemma}

{\sc Proof.} For $a>1$ let $f_a:{\bbR} \to {\bbR}$ be the convex increasing function defined by
$f_a(t)=t$ if $t<0$, $f_a(t)=at$ if $t \ge 0$. If $u\in PSH(\bbC^n)$ then
$f_a\circ u\in PSH(\bbC^n)$, $\gamma_{f_a\circ u}=a\gamma_u$, $\nu(f_a\circ u,z)=\nu(u,z)$ for any point $z$ and finally $\omega(S,f_a\circ u) = \omega(S,u)/a$. \cqfd\\

For technical reasons, we need to introduce some other constants. Let $l$ be a positive real, we define
$$\Omega_{psh}(S,l) :=\inf \{\gamma_u : u \in PSH(\bbC^n), \, \nu(u,p) \ge l \mbox{ for any } p \in S\}.$$
It is easy to see that  $\Omega_{psh}(S,l)/l$ doesn't depend on $l$. Then we just introduce the common value $\Omega_{psh}(S,1)$, that we denote by $\Omega_{psh}(S)$. 
In the same way we have 
$$\Omega^+_{psh}(S) :=\inf \{\gamma_u : u \in PSH(\bbC^n)\cap L^\infty_{loc}(\bbC^n \setminus S), \, \nu(u,p) \ge 1 \mbox{ for any } p \in S\}.$$
This invariant was already introduced in \cite{C-N}. All these invariants can be compared.

\begin{Lemma} \label{lm1.2}
{\it (i)} The sets $\{\gamma_u : u \in PSH(\bbC^n), \, \nu(u,p) \ge 1 \mbox{ for any } p \in S\}$ and $\{\gamma_u : u \in PSH(\bbC^n)\cap L^\infty_{loc}(\bbC^n \setminus S), \, \nu(u,p) \ge 1
\mbox{ for any } p \in S\}$ are connected and since the function $u(z)=\sum_{p \in S} \ln | z-p | + | z |$ satisfies $\gamma_u=+\infty$, we have  
$$(\Omega_{psh}(S),+\infty] \subset \{\gamma_u : u \in PSH(\bbC^n), \, \nu(u,p) \ge 1 \,\, \forall p \in S\} \subset [\Omega_{psh}(S),+\infty]$$ 
and 
$$(\Omega_{psh}^+(S),+\infty] \subset \{\gamma_u : u \in PSH(\bbC^n)\cap L^\infty_{loc}(\bbC^n \setminus S), \, \nu(u,p) \ge 1 \,\, \forall p \in S\} \subset [\Omega_{psh}^+(S),+\infty].$$
{\it (ii)} We have $\Omega_{psh}(S) \le \Omega(S)$ and $\Omega_{psh}(S) \le \Omega_{psh}^+(S)$.\\
{\it (iii)} We have $\omega(S) \ge \frac{| S |}{\Omega(S)}$, $\omega_{psh}(S) \ge \frac{| S |}{\Omega_{psh}(S)}$ and $\omega_{psh}^+(S) \ge \frac{| S |}{\Omega_{psh}^+(S)}.$\\
{\it (iii)}  If $S'\subseteq S\subset{\bbC}^n$ then
$$ \omega_{psh}^+(S') \le \omega_{psh}^+(S) \le \omega_{psh}^+(S') + \omega_{psh}^+(S\setminus S').$$
\end{Lemma}

The proof of this lemma is similar to the proof of Proposition $3.1$ in \cite{C-N}.\\
In what follows, we will prove that $\omega^+_{psh}(S) \le \Omega_{psh}(S) \le \Omega(S)$, in the $2$-dimensional case. So we are interested in finding lower bounds for $\omega_{psh}^+(S)$, 
and then for $\Omega(S)$. It is relevant for Conjecture $({\cal A}_1)$. On the other hand we are interested in finding upper bounds for $\omega_{psh}(S)$, 
and then for $\omega(S)$. It is relevant for Conjecture $({\cal A}_2)$, which is equivalent to Conjecture $({\cal A}_1)$.

\subsection{Comparisons of these invariants}

\subsubsection{A very simple situation in one complex variable}

The situation is very simple for subharmonic functions in $\bbC$. Let $v$ be an entire subharmonic (sh) function such that $\nu(v,p_j) = \alpha_j$ (positive real number) for $1\le j \le r$. 
Let $u$ be the sh function defined by $u(z)=\sum_{j=1}^r \alpha_j \ln | z -p_j |$.
Then $w:=v-u$ is a sh function in $\bbC$. 
$$0\le \gamma_w=\limsup_{ | z| \to +\infty} \frac{v(z)-u(z)}{\ln | z |}=
\limsup_{ | z| \to +\infty} \frac{v(z)}{\ln | z |} - \lim_{ | z| \to +\infty} \frac{u(z)}{\ln | z |}= \gamma_v - \gamma_u,$$ 
then $\gamma_v \ge \gamma_u = \sum_j \alpha_j $. $\omega(S,u)=1$ and $\omega(S,v)\le 1$. Consequently, $\omega_{sh}(S)= \omega_{sh}^+(S)=1$.\\
For polynomials, it is also very simple. Indeed, if $P$ is a polynomial such that $ord(P,p_j)=\alpha_j$ (positive integer) for $1\le j \le r$, then $P$ can be devided by 
$\prod_{j=1}^r (z-p_j)^{\alpha_j}$. There
exists a polynomial $Q$ such that $P(z)=\prod_{j=1}^r (z-p_j)^{\alpha_j} Q(z)$ and $deg(P) \ge deg(Q)$. Consequently, $\omega(S)=1$.\\
It is easy to see that we have also in this situation, $\Omega(S)=\Omega_{sh}(S) = \Omega_{sh}^+(S)=| S |$.

\subsubsection{In $\bbC^n$}

Here is a comparison principle which relates Lelong numbers at the points of $S$ of two psh functions with their logarithmic growth at infinity.
The proof of this result is simple and is similar to the proof of Theorem $3.4$ in $\cite{C-N}$ (or Proposition $2.1$ in \cite{Com1}).

\begin{Theorem} \label{thm3.2}
Let $S \subset \bbC^n$ be a finite set. Let $v$ be a psh function in $\bbC^n$. 
Let $u \in PSH(\bbC^n) \cap L^\infty_{loc}(\bbC^n \setminus S)$. Then
$$\sum_{p \in S} \nu(u,p)^{n-1} \nu(v,p) \le \gamma_u^{n-1} \gamma_v.$$
\end{Theorem}

\noindent {Proof.} Let $u \in PSH(\bbC^n) \cap L^\infty_{loc}(\bbC^n \setminus S)$ be such that $\gamma_u >0$. Let us suppose that $u$ has singularities at any point of $S$. \\
It follows from a comparison theorem for Lelong numbers with weights due to Demailly [Dem93] that for any $p \in S$
$$\lim_{r \to 0} \int_{B(p,r)} (dd^c u)^{n-1} \wedge dd^c v = \nu_u(dd^c v,p) \ge \nu(u,p)^{n-1} \nu(v,p)$$ 
and
$$\int_{\bbC^n} (dd^c u)^{n-1} \wedge dd^c v \ge \sum_{p \in S} \nu(u,p)^{n-1} \nu(v,p).$$
In addition,
$$\int_{\bbC^n} (dd^c u)^{n-1} \wedge dd^c v \le \gamma_u^{n-1} \int_{\bbC^n} (dd^c \log | z |)^{n-1} \wedge dd^c v \le \gamma_u^{n-1} \gamma_v.$$ 
and the proof is completed.\cqfd\\

\begin{Corollary} \label{cor3.2}
For any finite set $S \subset \bbC^n$, we have 
\begin{equation} \label{eq10}
\omega_{psh}^+(S)^{n-1} \le \Omega_{psh}(S) | S|^{n-2} \le \Omega(S) | S|^{n-2} \le | S|^{(n-1)^2/n},
\end{equation}
\begin{equation} \label{eq11}
\frac{| S |}{\Omega(S)} \le \omega(S) \le \omega_{psh}(S) \le \Omega_{psh}^+(S)^{n-1} \le | S |^{n-1}.
\end{equation}
In particular, for $n=2$, $\omega^+_{psh}(S) \le \Omega_{psh}(S)$ and $\omega_{psh}(S) \le \Omega_{psh}^+(S).$
\end{Corollary}

The inequalities $\omega(S) \le \Omega_{psh}^+(S)^{n-1} \le | S |^{n-1}$ are already proved in $\cite{C-N}$.\\ 

\noindent {\it Proof.} {\it (i)} Let $v$ be a psh function in $\bbC^n$ such that $\nu(v,p) \ge 1$ for any $p \in S$. Let $u$ be a function in $PSH(\bbC^n)\cap L^\infty_{loc}(\bbC^n \setminus S)$. According to theorem \ref{thm3.2}, we have 
$\displaystyle  \frac{\sum_{p \in S} \nu(u,p)^{n-1}}{\gamma_u^{n-1}} \le  \gamma_v.$ If we take the infimum in the right hand side of this inequality for any such $v$, then we obtain $\displaystyle \frac{\sum_{p \in S} \nu(u,p)^{n-1}}{\gamma_u^{n-1}}  \le \Omega_{psh}(S).$ By using H\"{o}lder inequality, we deduce $\displaystyle \frac{\sum_{p \in S} \nu(u,p)}{\gamma_u}  \le \Omega_{psh}(S)^{1/(n-1)} | S |^{(n-2)/(n-1)}$. Finally by taking the supremum in the left hand side of this inequality for any $u \in PSH(\bbC^n)\cap L^\infty_{loc}(\bbC^n \setminus S)$, we obtain
$$\omega_{psh}^+(S)^{n-1} \le \Omega_{psh}(S) | S|^{n-2}.$$
{\it (ii)} Let $v$ be a psh function in $\bbC^n$ and $u$ be a function in $PSH(\bbC^n)\cap L^\infty_{loc}(\bbC^n \setminus S)$ such that $\nu(u,p) \ge 1$ for any $p \in S$. According to theorem \ref{thm3.2}, we have $\displaystyle  \sum_{p \in S} \nu(v,p)/ \gamma_v \le  \gamma_u^{n-1}.$ We take the supremum in the left side of this inequality for any such $v$ and the infimum in the right side of this inequality for any such previous $u$ and we obtain 
$$\omega_{psh}(S) \le \Omega_{psh}^+(S)^{n-1}.$$ 
{\it (iii)} By definition of $\Omega_{psh}^+(S)$, with the psh function $u(z)=\sum_{p\in S} \log | z-p|$ we obtain 
$$\Omega_{psh}^+(S) \le | S |.$$
The others inequalities are obvious. This achieves the proof.\cqfd

\subsection{Proof of Theorem \ref{thm10.1}}




We already know that conjectures $({\cal P}_1)$ and $({\cal A}_1)$ are equivalent, according to Theorem \ref{thm2.1}. Indeed, with item $(iii)$, we have : $({\cal A}_1)$ implies $({\cal P}_1)$ and with item $(i)$, we have : $({\cal P}_1)$ implies $({\cal A}_1)$.\\
In addition it is well know that conjectures $({\cal A}_1)$ and $({\cal A}_2)$ are equivalent (\cite{Nag1}).\\
{\bf 1)} From conjecture $({\cal P}_3)$ we deduce conjectures $({\cal A}_1)$, $({\cal P}_2)$ and $({\cal A}_2)$.\\
Let us suppose that conjecture $({\cal P}_3)$ is satisfied: for any $\epsilon >0$, there exists an entire psh function $v$ in 
$L^{\infty}_{loc}(\bbC^n \setminus S)$, such that $\nu(v,p) \ge 1$ for any $p \in S$ and such that 
$\gamma_v \le (1+\epsilon)| S|^{1/n}$. Then $\omega(S,v)  \ge \frac{| S |}{(1+\epsilon)| S |^{1/n}}=\frac{| S |^{1-1/n}}{1+\epsilon}$. Consequently, 
\begin{equation} \label{eq12}
\omega_{psh}(S) \ge \omega_{psh}^+(S) \ge | S |^{1-1/n} 
\end{equation}
and 
\begin{equation} \label{eq13}
\Omega_{psh}(S) \le \Omega^+_{psh}(S) \le | S |^{1/n}.
\end{equation}
On the other hand, we have the chain $(\ref{eq10})$ of inequalities:
$$\omega_{psh}^+(S)^{n-1} \le \Omega_{psh}(S) | S |^{n-2} \le \Omega(S) | S |^{n-2} \le | S |^{(n-1)^2/n}.$$
Then we deduce that $\omega_{psh}^+(S)^{n-1} = \Omega_{psh}(S) | S |^{n-2} = \Omega(S) | S |^{n-2} =  | S |^{(n-1)^2/n},$
$$\omega_{psh}^+(S) = | S |^{(n-1)/n} \mbox{ and } \Omega_{psh}(S) =  \Omega(S) = | S |^{1/n}.$$
Conjecture  $({\cal A}_1)$  is solved.\\
According to $(\ref{eq13})$, we deduce that $\Omega_{psh}^+(S) = | S |^{1/n}.$ According to $(\ref{eq11})$, we obtain $\omega_{psh}(S) \le \Omega_{psh}^+(S)^{n-1} =| S |^{1-1/n}$. And according to inequalities $(\ref{eq12})$, we deduce that 
$$\omega_{psh}(S) = | S |^{1-1/n}.$$ 
Conjecture $({\cal P}_2)$ is solved.\\
Since $\omega(S) \le \omega_{psh}(S) = | S |^{1-1/n}$ and $\omega(S) \ge \frac{| S |}{\Omega(S)}$, conjecture $({\cal A}_2)$ is solved.\\
{\bf 2)} From conjecture $({\cal P}_2)$ we deduce conjectures $({\cal A}_1)$ and $({\cal A}_2)$.\\
Let us suppose that conjecture $({\cal P}_2)$ is satisfied: $\omega_{psh}(S) = \omega_{psh}^{+}(S)$. 
According to the two chains of inequalities $(\ref{eq10})$ and $(\ref{eq11})$, we obtain some identities.\\
$$\frac{| S |^{n-1}}{\Omega(S)^{n-1}} \le \omega(S)^{n-1} \le \omega_{psh}(S)^{n-1} \le \Omega_{psh}(S) | S |^{n-2} \le \Omega(S) | S |^{n-2} 
\le | S |^{(n-1)^2/n}.$$
$$\Arrowvert \hspace{4.5cm}$$
$$\omega_{psh}^{+}(S)^{n-1} \hspace{4.3cm} $$
We deduce in particular that $\displaystyle \frac{| S |^{n-1}}{\Omega(S)^{n-1}} \le \Omega(S) | S |^{n-2}$, which is equivalent to 
$\Omega(S) \ge | S |^{1/n}$. And since we always have $\Omega(S) \le | S |^{1/n}$, Conjecture $({\cal A}_1)$ is finally solved and all previous inequalities 
are equalities: 
$$| S |^{(n-1)^2/n} = \omega(S)^{n-1} = \omega_{psh}^{n-1}(S) = \omega_{psh}^{+}(S)^{n-1} = \Omega_{psh}(S) | S |^{n-2}.$$
We deduce that $\omega(S) = \omega_{psh}(S) = \omega_{psh}^{+}(S) = | S | ^{(n-1)/n}$ and $\Omega_{psh}(S)=| S | ^{1/n}$. In particular conjecture $({\cal A}_2)$ is solved.\\
{\bf 3)} Now let us prove that conjecture $({\cal A}_1)$ implies conjecture $({\cal P}_3)$. Since conjectures $({\cal A}_1)$ and $({\cal P}_1)$ are equivalent, $(g_1(tS,.))_{t\in \bbC^*}$ converges locally uniformly in $\overline{B(O,1)} \setminus \{O\}$ to $| S |^{1/n} \ln || . || $ : for any $\epsilon >0$ and any $0 < \varrho <1$, there exists $\epsilon_0 >0$ such that for any $| t | \le \epsilon_0$, we have 
$$(1+\epsilon)| S |^{1/n} \ln || z ||  \le g_1(tS,z) \le (1-\epsilon)| S |^{1/n} \ln || z ||  \,\,\,\mbox{ in } \overline{B(O,1)} \setminus B(O,\varrho).$$
Consequently, the following psh and continuous function $v$ is well defined in $\bbC^n$  
$$v(z)= \left\{ \begin{array}{lr}
g_1(tS,z), & z \in B(O,\varrho)\\
\mbox{max}\{g_1(tS,z), (1+\epsilon)| S |^{1/n} \ln || z ||\}, & z \in \overline{B(O,1)} \setminus B(O,\varrho)\\
(1+\epsilon)| S |^{1/n} \ln || z ||, & z \in \bbC^n \setminus  B(O,1)        
\end{array} \right.$$
Let us denote by $w$ the entire psh function defined by : $w(z)=v(tz)$. $w$ is locally bounded in $\bbC^n \setminus S$. $\nu(w,p)=1$ for any $p \in S$ and $\gamma_w = (1+\epsilon)| S |^{1/n}$, 
as required in conjecture $({\cal P}_3)$.
The proof is complete.\cqfd\\

\bibliographystyle{alpha}
\bibliography{NagataSHGH2}


\end{document}